\documentclass[12pt]{amsart}
\usepackage{graphicx}
\usepackage{amssymb}
\usepackage{epstopdf}
 \usepackage{latexsym}
\usepackage{amsmath}
\usepackage{amsfonts} 
\usepackage{amsthm}

\newtheorem{Proposition}{Proposition}
\newtheorem{Theorem}[Proposition]{Theorem}
\newtheorem{Lemma}[Proposition]{Lemma}

\newtheorem{Remark}{Remark}

    \def\sqr#1#2{{\vcenter{\vbox{\hrule height .#2pt
                             \hbox{\vrule width .#2pt height#1pt \kern#1pt
                                   \vrule width .#2pt}
                             \hrule height .#2pt}}}}

     \def\CC{\mathbb{C}}
    
    \def\NN{\mathbb{N}}

    \def\ZZ{\mathbb{Z}}

 \def\bfu{\mathbf{u}}
 \def\bfw{\mathbf{w}}
 \def\bff{\mathbf{f}}
  \def\bfg{\mathbf{g}}
 
 \def\bfy{\mathbf{y}}
  \def\bfh{\mathbf{h}}

 \def\bfq{\mathbf{q}}
   \def\bfm{\mathbf{m}}   
      
\def\bflam{\boldsymbol{\lambda}}
\def\bfphi{\boldsymbol{\phi}}
\def\bfpsi{\boldsymbol{\psi}}

\def\be{\begin{equation}}
\def\ee{\end{equation}}

\begin{document} 

\title[Differential systems with Fuchsian linear part]{Differential systems with Fuchsian linear part: correction and linearization, normal forms and multiple orthogonal polynomials}
\author{Rodica D. Costin}

\maketitle

\begin{abstract}
Differential systems with a Fuchsian linear part are studied in regions including all the singularities in the complex plane of these equations. Such systems are not necessarily analytically equivalent to their linear part (they are not linearizable) and obstructions are found as a unique nonlinear correction after which the system becomes formally linearizable. 

More generally, normal forms are found.

The corrections and the normal forms are found constructively. Expansions in multiple orthogonal polynomials and their generalization to matrix-valued polynomials are instrumental to these constructions.
\end{abstract}


\section{Introduction}

\subsection{The problem}

A first order linear system 
\be\label{eqw}
\frac{d\mathbf{w}}{dx}=A(x)\mathbf{w},\ \ \ {\mbox{with}}\ \bfw\in\CC^d,\ A(x)\in\mathcal{M}_{d}(\CC)
\ee
is called Fuchsian if all its singularities in 
$\overline{\CC}$ are regular. For such systems $A(x)$ has the form
\be\label{formA}
A(x)=\sum_{j=0}^{S+1}\frac{1}{x-p_j}A_j,\ \ \ {\mbox{with}}\ \bfw\in\CC^d,\ A_j\in\mathcal{M}_{d}(\CC)
\ee
and then (\ref{eqw}) is singular at $p_0,\ldots ,p_{S+1},\infty$, all of them regular singular points.\footnote{See \S\ref{sleq},\ref{App1s} for a classification of singularities.}$^,$\footnote{The peculiar choice of the limits of summation in (\ref{formA}) is for the convenience of the present paper.}  This linear context appears in a wide range of problems of mathematics and physics, and has been the topic of extensive studies.

The present paper studies nonlinear systems having a Fuchsian linear part:
\be\label{equ}
\frac{d\mathbf{u}}{dx}=A(x)\,\mathbf{u}+\frac{1}{Q(x)}\mathbf{f}(x,\mathbf{u}),\ \ \ \ \bfu\in\CC^d
\ee
with $A(x)$ given by (\ref{formA}). The term $\frac{1}{Q}\bff$ collects the nonlinear terms in $\bfu$ and $Q$ denotes the product
\be\label{notMQ}
Q(x)=\prod_{j=0}^{S+1}(x-p_j)
\ee
The question addressed is: under which conditions systems (\ref{equ}) are analytically equivalent to their linear part (\ref{eqw}) on a domain in $x$ large enough to encompass all the singularities $p_0,\ldots,p_{S+1}$?

This is a question of classification of differential equations in singular regions.
It is well known that in a neighborhood of a regular point, all equations are analytic equivalent. But this is no longer true on a domain including singularities.

\subsection{Motivation}

The question of linearization and, more generally, of equivalence, is a fundamental problem in the theory of differential equations.

 Besides the clear theoretical interest, linearization and equivalence are used as instrumental methods in control theory, and in devising algorithms for numerical and symbolic calculations (see, to cite just a few authors, \cite{De la Cru}, \cite{Doubrov}, \cite{Dridi}, \cite{Fels2}, \cite{Fels}). Often the methods used are developed from the method of equivalence introduced by Cartan \cite{Cartan} to decide whether two differential structures can be mapped one to another by a transformation taken in a given pseudogroup
\cite{Gardner}, \cite{Olver}. In the case of differential equations, the method was used for regular systems and extended for a neighborhood of one singular manifold \cite{Sternberg}, \cite{Hermann}.

There are many classes of problems that reduce to equations (\ref{formA}), (\ref{equ}). Vector fields with an eigenvalue equal to $1$ at a singular point can be reduced to the study of a Fuchsian system near one singularity by eliminating time.  Similarly, vector fields in regions containing two, or more singular points are reducible to (\ref{equ}).
More generally, Hamiltonian systems with polynomial potentials near particular, periodic, or doubly-periodic solutions can be reduced to (\ref{equ}) \cite{MAA}, \cite{Thesis}.

In dimension one nonlinear equations near a singular point (regular, or irregular) were studied in detail by Martinet and Ramis. Autonomous nonlinear equations near a singular point have been also studied, have generated deep results, and are relatively well understood. 

The present paper considers systems in regions encompassing two, or more, singularities.
The study also reveals a close connection between three concepts: linearizability, integrability and multiple orthogonal polynomials, and leads to generalizations of the latter.

The result of Theorem\,\ref{Th_obstr} has an interesting {{similarity with needed corrections found in other problems:}} \'Ecalle and Vallet showed that resonant systems are linearizable after appropriate correction \cite{Ecalle-Vallet}; also Gallavotti showed that there exists appropriate corrections of Hamiltonian systems so that the new system is integrable  \cite{Gallavotti}, convergence being proved later by Eliasson \cite{Eliasson}. This suggests the possible existence of an underlying general structure.

\subsection{Notations used}\

Denote $\NN=\{0,1,2,\ldots\}$.

For $\bfu,\,\bflam\in\CC^d$ and $\bfm\in\NN^d$, denote
$$\begin{array}{l}\bfu^\bfm=u_1^{m_1}u_2^{m_2}\ldots u_d^{m_d},\\|\bfm|=m_1+m_2+\ldots +m_d,\\\mathbf{m}\cdot {\bflam}=m_1\lambda_1+\ldots+m_d\lambda_d
\end{array}$$
$\mathcal{M}_{d}(\CC)$ denotes the $d\times d$ matrices with complex entries.\newline
$\CC^N[x]$ denotes the polynomials in $x$ with coefficients in $\CC^N$.\newline
$\mathcal{M}_N[x]$ denotes the polynomials in $x$ with coefficients in $\mathcal{M}_N(\CC)$.\newline
${\rm{e}}_i$ denotes the unit vector in $\CC^d$ with coordinate $i$ equal to $1$.

\subsection{Prior results}

\subsubsection{Systems near one singularity}\label{sysoneregsing}

If an equation (\ref{formA})-(\ref{equ}) is studied for $x$ in a domain small enough to include only one of the singularities $p_j$, and placing this singularity conventionally at $x=0$, such a system has the form
\be\label{non_one_sing}
\frac{d\mathbf{u}}{dx}=\frac{1}{x}L(x)\, \mathbf{u}+\frac{1}{x}\,\bff(x,\bfu),\ \ \ \bfu\in\CC^d
\ee
with $L(x)$ a matrix analytic at $0$ and $\bff(x,\bfu)$ analytic at $x=0$, with a zero of order two at $\bfu=0$. 

If the eigenvalues of $L(0)$ are linearly nonresonant, in the sense of (\ref{lnr}), then after an analytic change of coordinates the matrix $L(x)$ can be assumed constant $L(x)\equiv L(0)= L$ (see \S\ref{App1s}).

{Theorem}\,\ref{Th1s} shows that analytic linearization near one singularity holds for generic systems \cite{Nonln}:

\begin{Theorem}\label{Th1s}
Consider an equation
\be\label{my_th}
\frac{d{\mathbf{u}}}{dx}=\frac{1}{x}L\, {\mathbf{u}}+\frac{1}{x}\,\bff(x,\mathbf{u})
\ee
where $\bff(x,\bfu)$ is analytic for $|\bfu|<R$ ($R>0$) and for $x\in D$, where $D$ is a disk centered at $0$ (or an annulus centered at $0$).

It is assumed that $L$ is a diagonalizable matrix with  eigenvalues $\lambda_1,...,\lambda_d$ satisfying the following diophantine condition: there exist constants $C,\nu>0$ so that 
\be\label{Dio_cond}
\Big| k+ \bfm\cdot\boldsymbol{\lambda} -\lambda_i\Big|>C\left(|\bfm|+|k|\right)^{-\nu}
\ee
for all $i\in\{1,...,d\}$, $\mathbf{m}\in{\bf{N}}^d$, $|\mathbf{m}|\geq 2$, and for $ k\in\NN$ (respectively, $ k\in\ZZ$).

Then the system (\ref{my_th}) is analytically equivalent to its
linear part 
\be\label{lin1s}
\frac{d\mathbf{w}}{dx}=\frac{1}{x}L\, \mathbf{w}
\ee
 for $x\in D'\subset D$ and $|\bfw|<R'<R$.

Furhermore, $D'$ can be made arbitrarily close to $D$ if $R'$ is small enough.

Also, the analytic equivalence map is {unique} if no $\lambda_i$ is an integer.

\end{Theorem}

As a clear consequence of {Theorem}\,\ref{Th1s} the study of the local analytic properties of systems (\ref{my_th}) reduces to the study of the linear ones (\ref{lin1s}).

One application is to the study of integrability. Linear equations (\ref{lin1s}) have a full number of linearly independent first integrals which are analytic for small $x$ and $\bfw$, except for the singular manifolds $x=0$ and/or $\bfw=0$, where (generically) these first integrals have accumulation of poles: they are {\em{not meromorphic}} \cite{Nonln}. Then the same holds for the nonlinear systems (\ref{my_th}). 

A posteriori, it is clear that one should expect that first integrals, when they exist, are usually not meromorphic. On the other hand, when looking for first integrals in concrete problems, meromorphicity is often assumed. Such assumptions simplify the study, since under (some form of)  regularity assumptions the integrals are guaranteed to have specific types of expansions (e.g. Poincar\'e series, Laurent series) allowing for a study by successive approximations.

\subsubsection{Systems in regions with two regular singularities}

As seen in \S\ref{sysoneregsing}, near one singularity generic systems (\ref{formA})-(\ref{equ}) are linearizable, and integrable. Let us consider a larger domain of $x$, so that it includes two singularities. More precisely, 
consider (\ref{equ}) for $S=0$:
\be\label{equ2}
\frac{d\mathbf{u}}{dx}=\left(\frac{1}{x-p_0}A_0+\frac{1}{x-p_1}A_1\right)\mathbf{u}+\frac{1}{(x-p_0) (x-p_1)}\mathbf{f}(x,\mathbf{u})
\ee
 Of course, the position of the two singularities $p_{0},p_{1}$ can be arbitrarily placed in the complex plane using a linear change of variables.

The situation changes dramatically from the study near one singularity (\ref{my_th}): equations (\ref{equ2}) are not necessarily linearizable for $x$ in a domain $D$ containing both singularities $p_0$ and $p_1$. 

This fact can be understood by applying Theorem\,\ref{Th1s} near singularity $p_0$, and then, near singularity $p_1$. Then (generically\footnote{Existence of an analytic lienarization map for $x$ close to $p_j$ follows from Theorem\,\ref{Th1s} if $A_j$ is diagonalizable, satisfies the Diophantine condition (\ref{Dio_cond}), and is linearly nonresonant, \S\ref{nl1fuchs}.})
there exists a {unique} linearization map analytic at $x=p_0$, and
there exists a {unique} linearization map analytic at $x=p_1$. But the two maps do not necessarily coincide, as the map analytic at $x=p_j$ is usually branched
 at $x=p_{1-j}$ ($j=0,1$).

This fact has important consequences, as it turns out that linearizability and integrability are intimately connected \cite{RDC-MDK} : 

\begin{Theorem}

Consider an equation (\ref{equ2}) in dimension 1:

\be\label{eq1d2s}
\frac{du}{dx}=\left(\frac{a_0}{x-1}+\frac{a_1}{x+1}\right)u+\frac{f(x,u)}{x^2-1}
\ee
 
 If equation  (\ref{eq1d2s}) is {\em{not}} analytically linearizable then no single-valued integrals exists on domains encircling $\pm1$ for generic\footnote{Precise conditions are given in \cite{RDC-MDK}.} $a_0,a_1$.

\end{Theorem}

It is then important to determine which equations are linearizable, and which are
not, and to find the equivalence classes with respect to analytic equivalence in this semi-local context.

This problem was studied in \cite{RDC_classification} for $S=0$; the results follow from Theorem\,\ref{Th_obstr} proved here.

\section{Main results}

\subsection{Setting}\label{setting}
Let $S\geq 0$, let $p_0,p_1,\ldots ,p_{S+1}$ be {\em{distinct}} points in $\CC$ and let $A_0,A_1,\ldots,A_{S+1}$ be $d\times d$ matrices. 

Consider the differential system (\ref{equ}) with (\ref{formA}), (\ref{notMQ}).

Let $D$ be a simply connected domain in the complex plane containing the singularities $p_0,\ldots,p_{S+1}$.

It is assumed that $\bff(x,\cdot)$ has a zero of order two at $\bfu=0$, 
and is assumed analytic for $x\in D$ and small $|\bfu|$. 

The denominator $Q$ of the nonlinear part $\frac{1}{Q}\bff$ of (\ref{equ}) simply shows that this nonlinear part may have at most poles of order one at $p_0,\ldots,p_{S+1}$.

Denote
\be\label{Ainfi}
A_\infty=\sum_{j=0}^{S+1}\, A_j
\ee
Note that $A(x)\sim A_\infty x^{-1}$ as $x\to\infty$.

{\bf{Definition:}} A matrix $M\in\mathcal{M}_{d}(\CC)$ is called nonlinearly nonresonant if its eigenvalues $\lambda_1,\ldots,\lambda_d$ satisfy
\be\label{condA}
k+\boldsymbol{\lambda}\cdot\bfm-\lambda_i\ne 0,\ {\mbox{for\ all\ }} k\in\NN,\  \bfm\in\NN^d,\ |\bfm|\geq 2,\ i=1,\ldots,d
\ee

{\bf{Assumption:}}
It is assumed that each matrix $A_0,A_1,\ldots,A_{S+1},A_\infty$ is nonlinearly nonresonant.

\subsection{Main results}

\begin{Theorem}\label{Th_obstr}

Consider an equation (\ref{formA})-(\ref{equ}) under the notations and assumptions of \S\ref{setting}.

Then there exists a {\em{unique}} correction of the nonlinear part as a formal series
\be\label{serphi}
\boldsymbol{\phi}(x,\bfu)=\sum_{\bfm\in\NN^d,\ |\bfm|\geq 2}\, \boldsymbol{\phi}_\bfm(x)\bfu^\bfm
\ee
where $\boldsymbol{\phi}_\bfm(x)$ are {\em{polynomials in $x$ of degree at most $S$}}, such that the corrected system
\be\label{eqSc}
\frac{d\mathbf{u}}{dx}=A(x)\,\mathbf{u}+\frac{1}{Q(x)}\left[\mathbf{f}(x,\mathbf{u})-\boldsymbol{\phi}(x,\bfu)\right],\ \ \ \ \bfu\in\CC^d
\ee
is formally linearizable through a change of variable as a series
\be\label{serh}
\bfu=\bfw+\bfh(x,\bfw)=\bfw+\sum_{|\bfm|\geq 2}\, \bfh_\bfm(x)\bfw^\bfm
\ee
with $\bfh_\bfm(x)$ analytic on $D$.

\end{Theorem}

{\bf{Remarks.}}

{\bf{1.}} Equation (\ref{equ}) is formally linearizable if and only if $\boldsymbol{\phi}(x,\bfu)\equiv 0$. Therefore this unique correction gathers the obstructions to linearizability.

{\bf{2.}} Obviously, if an equation is not formally linearizable, then it is not analytically linearizable either.

Since equations (\ref{equ}) are not necessarily linearizable, then they are not all equivalent either. Theorem\,\ref{Th_norm_form} provides the classification of these equations by specifying formal normal forms:

\begin{Theorem}\label{Th_norm_form}

Consider an equation (\ref{formA}),(\ref{equ}) under the notations and assumptions of \S\ref{setting}.

Then there exists a {\em{unique}}  formal series
\be\label{serpsi}
\boldsymbol{\psi}(x,\bfu)=\sum_{\bfm\in\NN^d,\ |\bfm|\geq 2}\, \boldsymbol{\psi}_\bfm(x)\bfu^\bfm
\ee
where $\boldsymbol{\psi}_\bfm(x)$ are {\em{polynomials in $x$ of degree at most $S$}}, such that the system (\ref{equ}) is formally equivalent to
\be\label{norm_from}
\frac{d\mathbf{w}}{dx}=A(x)\,\mathbf{w}+\frac{1}{Q(x)}\boldsymbol{\psi}(x,\bfw),\ \ \ \ \bfw\in\CC^d
\ee
through a change of variable as a series of the form (\ref{serh})
with $\bfh_\bfm(x)$ analytic on $D$.

\end{Theorem}


{\bf{Comments about convergence.}} While $\bfphi_\bfm(x)$, $\bfpsi_\bfm(x)$ are polynomials, and 
$\bfh_\bfm(x)$ are analytic functions, the series in $\bfw$ (\ref{serphi}), (\ref{serh}), (\ref{serpsi}) are formal. 

Convergence of (\ref{serphi}) and of the linearization (\ref{serh}) were proved in \cite{Analytic} for $S=0$ under the supplementary assumptions that $A_0,A_1$ are simultaneously diagonalizable, and have eigenvalues with positive real parts. 

For convergence to hold in general it is reasonable to expect that the nonresonance condition (\ref{condA}) should be strengthened to a diophantine condition. This could be a Siegel type condition, like (\ref{Dio_cond}), or even a Bruno type condition; the latter was proved by Yoccoz to be optimal for the problem of iteration of holomorphic germs in one variable, and under this condition convergence holds in the linearization problem for vector fields at a resonant point, as proved by \'Ecalle and Vallet \cite{Ecalle-Vallet}.

\section{Proofs}

\subsection{{The structure of the proofs.}} 

The main steps in the proofs of Theorems\, \ref{Th_obstr} and\,\ref{Th_norm_form} are as follows. Using power series expansions in $\bfw$ it turns out that the terms $\bfh_\bfm(x)$ in (\ref{serh}) satisfy recursive systems of differential equations which are first order, linear nonhomogeneous. Remarkably, their linear part is Fuchsian. The burden of proof is carried by the study of these Fuchsian nonhomogeneous systems, and the structure of their solutions turns out to be rich and interesting.


\subsection{A fundamental lemma}\label{fdLem}

Consider a Fuchsian equation with a nonhomogeneous term:
\be\label{prototype}
\bfy'(x)+\, B(x)\, \bfy(x)=\, \frac{1}{Q(x)}\, \bfg(x),\ \ \ \ \ \ \ \ \ \bfy,\bfg\in\CC^N,\ B(x)\in\mathcal{M}_{N}(\CC)
\ee
where
\be\label{fuchsian}
 B(x)=\sum_{j=0}^{S+1}\frac{1}{x-p_j}B_j,\ \ \ \ \ \ Q(x)=\prod_{j=0}^{S+1}(x-p_j)
 \ee
Let
\be\label{not_gamma}
B_\infty=\sum_{j=0}^{S+1}B_j
\ee
It is assumed that
\be\label{kpBinvert}  
k+B_j\ \ {\mbox{are\ invertible\ for\ all\ }} k\in\NN\ {\mbox{and}}\ j=0,1,\ldots,S+1,\infty .
\ee

\begin{Lemma}\label{FundLem} 

Consider the equation (\ref{prototype})-(\ref{fuchsian}) under the assumption (\ref{kpBinvert}).

Let $D$ be a simply connected domain containing the singularities $p_0,\ldots,p_{S+1}$.

Then for any function $\bfg(x)$ analytic on $D$ there exists a unique $\boldsymbol{\phi}(x)\in\CC^N[x]$ polynomial of degree at most $S$ so that the corrected equation
 \be\label{corprototype}
\bfy'(x)+\, B(x)\, \bfy(x)=\, \frac{1}{Q(x)}\, \left[\bfg(x)-\boldsymbol{\phi}(x)\right]
\ee
has a solution $\bfy(x)$ which is analytic on $D$.

\end{Lemma}

The proof of Lemma\,\ref{FundLem} is the topic of \S\ref{PFL}. Before providing its details, which involve both analytic and algebraic arguments, and leads to generalizations of multiple orthogonal polynomials to matrix-valued ones, let us assume {Lemma}\,\ref{FundLem} true for the moment and show its usefulness by providing the proofs of Theorems\,\ref{Th_obstr} and\,\ref{Th_norm_form}.

\subsection{Proof of Theorem\,\ref{Th_obstr}}


The change of variables (\ref{serh}) provides a linearization of (\ref{eqSc}) if and only if $\bfh$ satisfies the nonlinear partial differential equation
\be\label{pdel1}
\partial_x\bfh+d_\bfw\bfh\, A\bfw=A\bfh+\frac{1}{Q(x)}\, \left[\bff(x,\bfw+\bfh)-\bfphi(x,\bfw+\bfh)\right]
\ee

Searching for solutions of (\ref{pdel1}) as power series in $\bfw$, denote by $\bfh_n$ the homogeneous part degree $n$ of $\bfh(x,\bfw)$; in the notation (\ref{serh}) we have
\be
\bfh_n(x,\bfw)=\sum_{|\mathbf{m}|=n}\bfh_\mathbf{m}(x)\bfw^\mathbf{m}
\ee
We use similar notations to denote by $\bff_n$ and $\bfphi_n$ the homogeneous part of degree $n$ of $\bff$, respectively $\bfphi$.

Equation (\ref{pdel1}) splits into blocks of systems of ordinary differential equations for $\{\bfh_\bfm\}_{|\bfm|=n}$, recursive in $n$:
\be\label{eqhn}
\partial_x\bfh_n+{\rm{d}}_\bfw\bfh_n\, A\bfw-A\bfh_n=\frac{1}{Q(x)}\, \mathbf{R}_n(x,\bfw),\ \ \ \ n\geq 2
\ee
where
\be\label{formRn}
 \mathbf{R}_n=\bff_n-\bfphi_n+\tilde{\mathbf{R}}_n
\ee
with $\tilde{\mathbf{R}}_n$ a
polynomial in $\bfphi_\mathbf{m}$, $\bfh_\bfm$, $\bff_\mathbf{m}$ with $|\mathbf{m}|<n$, and $\tilde{\mathbf{R}}_2=0$.

Each $\bfh_n$ and $\bfphi_n$ are to be determined from (\ref{eqhn})-(\ref{formRn}) inductively on $n$. This can be done as follows.

Note that, remarkably, for each $n$, (\ref{eqhn}) is a Fuchsian nonhomogeneous system.

To see this, denote by $\mathcal{P}_n$ the space of vector-valued polynomials in $\bfw\in\CC^d$, homogeneous degree $n$:
\be\label{defPn}
 \mathcal{P}_n\, =\, \left\{\ \bfq\ ;\ \bfq(\bfw)=\sum_{\bfm\in\NN^d, |\bfm|=n}\bfq_\bfm\bfw^\bfm,\ \bfq_\bfm\in\CC^d\ \right\}
 \ee
Note the canonical basis of the linear space $ \mathcal{P}_n$: 
\be\label{canbasP}
\mathbf{r}_{\bfm,i}=\bfw^\bfm{\rm{e}}_i,\ \ |\bfm|=n,\ i=1,\ldots, d
\ee
and its finite dimension $dim\,\mathcal{P}_n=N=N(d,n)$.\footnote{In fact $N=d(n+d-1)!/n!/(d-1)!$.}

Denote by $B(x)$ the linear operator on  $ \mathcal{P}_n$ taking $\bfh_n$ to 
$$B(x)\bfh_n={\rm{d}}_\bfw\bfh_n\, A\bfw-A\bfh_n$$
 We have
\be\label{BF}
B(x)=\sum_{j=0}^{S+1}\frac{1}{x-p_j}\, B_j\ \ \ \ \ {\mbox{where}}\ \ B_j\bfq={\rm{d}}_\bfw\bfq\, A_j\bfw-A_j\bfq
\ee
hence in the canonical base (\ref{canbasP}) $B(x)$ is given by a Fuchsian matrix: for each $n$ the system (\ref{eqhn}),(\ref{formRn}) has the form (\ref{corprototype}) (with $\bfg=\bff_n+\tilde{\mathbf{R}}_n $ and $\bfphi=\bfphi_n$). 

Note that $B_\infty={\rm{d}}_\bfw\bfq\, A_\infty\bfw-A_\infty\bfq$.

Theorem\,\ref{Th_obstr} now follows from Lemma\,\ref{FundLem} once we show that its assumption
(\ref{kpBinvert}) is satisfied. This is proved by the following result:

\begin{Lemma}\label{sollem} 

Let $M$ be a $d\times d$ matrix with eigenvalues $\lambda_1,\ldots,\lambda_d$ (not necessarily distinct). 

Let $J_M$ be the following linear operator on $\mathcal{P}_n$:
\be\label{defJ}
(J_M\bfq)(\bfw)=({\rm{d}}_\bfw\bfq )\, M \bfw-M\bfq
\ee

Then the eigenvalues of $J_M$ are $\boldsymbol{\lambda}\cdot\bfm-\lambda_i$, $i=1,\ldots,d$, $|\bfm|=n$.



\end{Lemma}

The proof of Lemma\,\ref{sollem} is found in \S\ref{proofLJ}. 

In the notation (\ref{BF}), (\ref{defJ}) we have $B_j=J_{A_j}$ for all $j=0,\ldots,S+1,\infty$. Therefore, in view of Lemma\,\ref{sollem}, (\ref{condA}) implies (\ref{kpBinvert}).\qed

\subsection{Proof of Theorem\,\ref{Th_norm_form}}


The change of variables (\ref{serh}) transforms (\ref{eqSc}) into (\ref{norm_from}) if and only if $\bfh$ satisfies the nonlinear partial differential equation
\be\label{pdel}
\partial_x\bfh+d_\bfw\bfh\, \left( A\bfw+\frac{1}{Q}\bfpsi(x,\bfw)\right)=A\bfh+\frac{1}{Q}\, \left[\, \bff(x,\bfw+\bfh)-\bfpsi(x,\bfw)\right]
\ee

Expanding in power series in $\bfw$ we obtain, as in the proof of Theorem\,\ref{Th_obstr}, that the homogeneous parts $\bfh_n$ ($n\geq 2$) of $\bfh$ satisfy systems of equations which can be solved recursively in $n$:

\be\label{eqhn2}
\partial_x\bfh_n+{\rm{d}}_\bfw\bfh_n\, A\bfw-A\bfh_n=\frac{1}{Q(x)}\, \mathbf{R}_n(x,\bfw),\ \ \ n\geq 2
\ee
where
\be\label{formRn2}
 \mathbf{R}_n=\bff_n-\bfpsi_n+\tilde{\mathbf{R}}_n
\ee
with $\tilde{\mathbf{R}}_n$ a
polynomial in $\bfpsi_\mathbf{m}$, $\bfh_\bfm$, ${\rm{d}}_\bfw\bfh_\bfm$, $\bff_\mathbf{m}$ with $|\mathbf{m}|<n$, and $\tilde{\mathbf{R}}_2=0$.

Each $\bfh_n$ and $\bfphi_n$ can now be determined from (\ref{eqhn2}) as in the proof of Theorem\,\ref{Th_obstr}.\qed

\section{Proof of Lemma\,\ref{FundLem}}\label{PFL}

The proof consists of several parts. Section \S\ref{Step0} shows that if the correction $\bfphi$ exists, then it is unique. 

Section \S\ref{Step_I} proves Lemma\,\ref{FundLem} for the case when $\bfg(x)$ is a polynomial. The proof is algebraic, and solutions are found as expansions in terms of multiple orthogonal polynomials and their generalizations introduced here.

Section \S\ref{Step_II} proves Lemma\,\ref{FundLem} for the case when the eigenvalues of all matrices $B_j$ have positive real parts. An analytic approach is used.

Section \S\ref{Step_III} shows that the general case can be reduced to \S\ref{Step_II}, thus completing the proof of Lemma\,\ref{FundLem}. The algebraic results of  \S\ref{Step_I} are used to construct this reduction.

\subsection{Existence of local analytic solutions.}\label{locsol}

The existence and uniqueness of an analytic at $p_j$ solution of  (\ref{prototype}) follows from the rest of \S\ref{PFL}. Due to its own intrinsic interest the statement and a simple proof are included here.

\begin{Lemma}\label{L7}
Under the assumption (\ref{kpBinvert})
equation (\ref{prototype})-(\ref{fuchsian}) has a unique solution $\bfy(x)$ which is analytic at the singular point $x=p_j$ (for $0\leq j\leq S+1$). 
\end{Lemma}

{\em{Proof of Lemma\,\ref{L7}.}}

Substituting a power series $\bfy(x)=\bfy_0+(x-p_j)\bfy_1+(x-p_j)^2\bfy_2+\ldots$ in (\ref{prototype}) one obtains a recursive system for the coefficients $\bfy_0,\bfy_1,\bfy_2,\ldots$:
$$\begin{array}{ll}
B_j\bfy_0=\bfg_0\\
\displaystyle{(1+B_j)\bfy_1=\bfg_1-\, \sum_{k;\, k\ne j}\, \frac{1}{p_j-p_k}B_k\bfy_0}\\
\vdots\\
(k+B_j)\bfy_k=\mathbf{r}_k(\bfy_0,\ldots,\bfy_{n-1},\bfg)\\
\vdots
\end{array}
$$
which has a unique solution by the assumption (\ref{kpBinvert}). 

The series converges since equation (\ref{prototype}) is linear nonhomogeneous.\qed

\subsection{Uniqueness of corrections}\label{Step0}

If there are two corrections $\bfphi=\bfphi_1$, and $\bfphi=\bfphi_2$,  so that equation
 (\ref{corprototype}) has the analytic on $D$ solutions $\bfy=\bfy_1$, respectively $\bfy=\bfy_2$, then $\bfy=\bfy_1-\bfy_2$ is an analytic solution of (\ref{prototype}) with $\bfg=\bfphi_2-\bfphi_1$ a polynomial degree at most $S$. Lemma\,\ref{UniLem} shows that this implies $\bfphi_1=\bfphi_2$ and $\bfy_1=\bfy_2$.
 
The following result is needed.
 
 \begin{Lemma}\label{pgr}
 
If $\bfg(x)$ is a polynomial then any solution of (\ref{prototype})-(\ref{fuchsian}) is polynomially bounded for $x\to\infty$. 
 
 \end{Lemma}
 
 {\em{Proof of Lemma\,\ref{pgr}.}} 
 
Denote by $W(x)$ the integrating factor of (\ref{prototype}): an invertible matrix satisfying
\be\label{defW}
\frac{d}{dx}W(x)=W(x)B(x)
\ee
Note that $W(x)$ is a fundamental matrix for the equation 
\be\label{eqv}
\mathbf{v}'=\mathbf{v} B(x)
\ee 
and that $W$ is the inverse of a fundamental matrix of the linear equation $\bfy'+B(x)\bfy=0$.

Since for (\ref{eqv}) the point at $\infty$ is Fuchsian the matrix $W$ can be chosen such that 
$W(x)=x^{\hat{B}_\infty}\Phi(x)$ where $\hat{B}_\infty$ is a constant matrix whose eigenvalues (not counting the multiplicity) are also eigenvalues of ${B_\infty}$, with $\Phi(x)$ analytic at $\infty$ and $\Phi(x)=I+o(1)$ as $x\to\infty$ (see \S\ref{App1s} for details). Therefore
\be\label{chooseW}
W(x)\, =\, x^{\hat{B}_\infty}\left(I+o(1)\right)\ \ \ {\mbox{for\ }}x\to\infty
\ee

The general solution of (\ref{prototype}) is
\be\label{gsnsolln}
W(x)^{-1}\, \boldsymbol{\xi}\, +\, W(x)^{-1}\int^x\, Q(t)^{-1}W(t)\bfg(t)\, dt
\ee
where the path of integration avoids the singularities $p_j$, the integral starts at an arbitrary, but fixed point, and $\boldsymbol{\xi}$ is a constant vector.

Denote $\bfg(x)=\bfg_mx^m+\ldots +\bfg_1x+\bfg_0$.

For $x\to\infty$ we have $W(x)\sim x^{\hat{B}_\infty}$, $Q(x)\sim x^{S+2}$ and $\bfg(x)\sim \bfg_mx^m$ therefore (\ref{gsnsolln}) is polynomially bounded for $x\to\infty$. \qed

\begin{Lemma}\label{UniLem} 

Consider equations (\ref{prototype})-(\ref{fuchsian}) under the assumption (\ref{kpBinvert}).

If $\bfg(x)$ is a polynomial of degree at most $S$ then equation 
(\ref{prototype}) has a solution which is analytic at all the singular points $p_0,\ldots , p_{S+1}$ only for $\bfg(x)\equiv 0$, and in this case the only analytic solution is $\bfy(x)\equiv 0$.

\end{Lemma}

{\em{Proof of Lemma\,\ref{UniLem}.}}

Assume that there exists $\bfy(x)$ a solution of (\ref{prototype}) which is analytic at all the singular points $p_0,\ldots , p_{S+1}$ of the equation. Since (\ref{prototype}) is a linear nonhomogeneous ordinary differential equation then $\bfy(x)$ can have no other singularities and therefore $\bfy(x)$ is entire. By Lemma\,\ref{pgr} then $\bfy(x)$ is a polynomial. 

If $\bfy(x)\equiv 0$ then (\ref{prototype}) implies $\bfg(x)\equiv 0$.

Assume, to obtain a contradiction, that $\bfy(x)\not\equiv 0$. Denote $\bfy(x)=\bfy_nx^n+\ldots+\bfy_0$ with $\bfy_n\ne 0$ and $\bfg(x)=\bfg_mx^m+\ldots +\bfg_1x+\bfg_0$.
Substituting these expansions in (\ref{prototype}), for $x\to\infty$ the left-hand-side has the dominant term $(n+{B_\infty} )\bfy_nx^{n-1}$ and the right hand-side has the dominant term $\bfg_mx^{m-S-2}$. Since $m\leq S$ then $m-S-2<n-1$ therefore we must have $(n+{B_\infty} )\bfy_n=0$, and by assumption (\ref{kpBinvert}) then $\bfy_n$=0, which is a contradiction. \qed

\subsection{Algebraic structure: proof of Lemma\,\ref{FundLem} for $\bfg(x)$ polynomials using  expansions in multiple orthogonal polynomials and their generalizations}\label{Step_I}

\

In this section it is assumed that $\bfg(x)$ are polynomials.
Lemma\,\ref{FundLem}, stating existence of corrections and of corresponding analytic solutions (which turn out to be polynomials in this case) is proved using expansions in multiple orthogonal polynomials and their generalizations, introduced here.

Let $B(x)$ be a Fuchsian matrix (\ref{fuchsian}), with the notation (\ref{not_gamma}). Condition 
(\ref{kpBinvert}) is not used until \S\ref{422}.

Let $W(x)$ be a fundamental matrix of (\ref{eqv}). 

For any $n=0,1,2,\ldots$ denote 
\be\label{notnmk}
n=(S+1)m+i\ {\mbox{ with\ }}m=\lfloor n/(S+1)\rfloor,{\mbox{and}}\ i=0,1,\ldots , S
\ee
A precise notation requires $m=m_n$ and $i=i_n$, but the subscript will be omitted to avoid a heavy notation. Anyhow, the one-to-one correspondence $n\to(m,i)$ implied by (\ref{notnmk}) is only used in \S\ref{Step_I}.

With the notation (\ref{notnmk}), define
\be\label{defPn}
P_{n}^{(W)}\,(x)\equiv P_{n}\,(x)=W(x)^{-1}\, \frac{d^m}{dx^m}\, \left[\, x^i\, Q(x)^m\, W(x)\,\right]
\ee

{\bf{Remarks.}}

{\bf{1.}} In dimension one ($N=1$) we have $B_j=b_j\in\CC$, and $W(x)$ has the form
\be\label{defW1d}
W(x)=\prod_{j=0}^{S+1}(x-p_j)^{b_j}
\ee

{\bf{1a.}} In particular, for $S=0$ only $i=0$ can occur in (\ref{notnmk}), and  (\ref{defPn}) is the Rodrigues formula defining the Jacobi polynomials $P^{(b_0,b_1)}_n(x)$ (up to a multiplicative constant).

{\bf{1b.}}   For $S=1$ formula (\ref{defPn}) defines the Jacobi-Angelesco multiple orthogonal polynomials introduced and studied in \cite{Kaliaguine} (see also \cite{Ismail}). 

{\bf{1c.}}  For higher $S$ the $P_n(x)$ defined by (\ref{defPn}) are a natural extension of the Jacobi-Angelesco polynomials. 

{\bf{2a.}} For higher dimensions ($N>1$) and $S=0$ formula (\ref{defPn}) was studied in \cite{GenPol} where it is shown that $P_n(x)$ are polynomials and that they form a complete system. In addition they satisfy a two step recurrence relation, and, in the commutative case, other beautiful properties of orthogonal polynomials also hold: they are eigenfunctions for a second order differential operator, and they are orthogonal with respect to the weight $W(x)$ (when integrals exist).

{\bf{2b.}} For $S=0$ the polynomials  occurring in the proof of Theorem\,\ref{Th_obstr} were studied in 
\cite{Gen_Jacobi}.

\

For the proofs in the present paper we need to show that formula (\ref{defPn}) defines polynomials indeed (done in Proposition\,\ref{prod_form}), that they form a complete system 
({Proposition}\,\ref{polyno}), and that (\ref{corprototype}) is solved in terms of them, \S\ref{422}.

\subsubsection{$P_n(x)$ are polynomials.}\label{poly}

{{Note}} first that $P_i(x)=x^i$ for $i=0,1,\ldots S$.

The following reformulation of Rodrigues-type formula  (\ref{defPn}) helps simplify the proofs considerably (the case $S=0$ appears in \cite{GenPol}).

We let $\ \ \displaystyle{\frac{d}{dx}=\partial_x}$.

\begin{Proposition}\label{prod_form}
Denote by $\mathcal{A}_k$ the following linear operators on $\mathcal{M}_N[x]$:
\be\label{defAk}
\mathcal{A}_k\, =\, k\, Q'(x)\, +Q(x)B(x) +\, Q(x)\, \partial_x
\ee

Then the functions $P_n$ defined by (\ref{defPn}) satisfy
\be\label{newPn}
 P_n(x)=\mathcal{A}_1\mathcal{A}_2\ldots\mathcal{A}_m\, x^i
 \ee
 for $n\geq S+1$. 
 
 As a consequence, $P_n(x)$ are polynomials in $x$ (matrix-valued).
\end{Proposition}

{\em{Proof of {Proposition}\,\ref{prod_form}.}} 

Calculating successively  the $m$ derivatives in (\ref{defPn}) we obtain, using (\ref{defW}) and 
(\ref{defAk}),
\begin{multline}
\frac{d}{dx}\left[\, x^iQ^m\, W\, \right]\, 
=\, mx^iQ^{m-1}Q'W+\, x^iQ^mW'\, +Q^mW\partial_xx^i\\
=\, Q^{m-1}W\left[\left( mQ'+\, QW^{-1}W'\right)\, x^i\, +Q\partial_xx^i\right]\, =Q^{m-1}W\mathcal{A}_m\, x^i
\end{multline}
then
\begin{multline}
\frac{d^2}{dx^2}\left[\, x^iQ^m\, W\, \right]\, =\, \frac{d}{dx}\left[\, Q^{m-1}\, W\, \mathcal{A}_m\, x^i\right]\, \\
=\, Q^{m-2}\, W\, \mathcal{A}_{m-1}\mathcal{A}_m\, x^i
\end{multline}
and so on.\qed

\begin{Proposition}\label{polyno}

(i) The polynomial $P_n$ defined by (\ref{defPn}) has degree at most $n$, and the coefficient of $x^n$ is 
\be\label{domcf}
C_n=\prod_{j=1}^m\left( j+n+B_\infty\right)\ \ \ {\mbox{for}}\ n\geq S+1
\ee

(ii) If $B_\infty$ satisfies
\be\label{wib}
k+B_\infty{\mbox{\ is invertibe\ for\ all\ }}k\geq S+2,\ k\in\NN
\ee
then (\ref{domcf}) is an invertible matrix, and therefore the degree of $P_n$ is exactly $n$.

As a consequence, the set $\{P_n(x)\}_{n\in\NN}$ is complete in the sense that any polynomial $\bff(x)\in\CC^N[x]$ can be written as an expansion in $P_n$'s:
\be\label{expg}
\bff(x)\, =\, \sum_{n=0}^{{\rm{deg}}\,\bff}\, P_n(x)\bff_n,\ \ \ \ {\mbox{for\ some\ }}\bff_n\in\CC^N
\ee
and the expansion is unique.

\end{Proposition}

{\em{Proof of {Proposition}\,\ref{polyno}.}}

{\em{(i)}} Retain only the dominant coefficients in (\ref{defAk}) using the fact that $Q(x)\sim x^{S+2}$ and $B(x)\sim  x^{-1}{B_\infty}$ as $x\to\infty$. We then have
$$\mathcal{A}_k\sim x^{S+1}\left[k(S+2)+B_\infty+x\partial_x\right]$$
Therefore
$$\mathcal{A}_mx^i\sim \left[ m(S+2)+B_\infty+i\right]x^{i+S+1}$$
then
\begin{multline}\nonumber
\mathcal{A}_{m-1}\mathcal{A}_mx^i\sim \\
\left[(m-1)(S+2)+B_\infty+i+S+1\right]\, \left[m(S+2)+B_\infty+i\right]x^{i+2(S+1)}
\end{multline}
and so on. Using (\ref{newPn}) therefore
$$P_n(x)\sim\prod_{j=1}^m\left[j(S+2)+B_\infty+i+(m-j)(S+1)\right]x^{i+m(S+1)}$$
which equals $C_nx^n$, proving (\ref{domcf}). 

{\em{(ii)}} Condition (\ref{wib}) clearly ensures that (\ref{domcf}) is an invertible matrix. Invertibility of the dominant coefficient of $P_n$ implies, using induction over the degree of $\bff$, the existence of the unique expansion (\ref{expg}).\qed


\subsubsection{Proof of Lemma\,\ref{FundLem} for $\bfg(x)$ polynomial}\label{422}

Assume (\ref{kpBinvert}) satisfied.

Denote $\tilde{\tilde{W}}=\frac{1}{Q}W$. Note that $\tilde{\tilde{W}}$ is a fundamental matrix of solutions for the Fuchsian equation $\tilde{\tilde{\mathbf{v}}}'=\tilde{\tilde{\mathbf{v}}}\tilde{\tilde{B}}$ where $\tilde{\tilde{B}}=B-Q'/Q$. Note that $\tilde{\tilde{B}}=\sum\frac{1}{x-p_j}\tilde{\tilde{B}}_j$ with $\tilde{\tilde{B}}_j=B_j-I$, and therefore
\be\label{ttbi}
\tilde{\tilde{B}}_\infty=B_\infty-(S+2)I
\ee

Consider the polynomials associated to $\tilde{\tilde{W}}$ by a formula (\ref{defPn}): $P^{(\tilde{\tilde{W}})}_n\equiv\tilde{\tilde{P}}_n$. Note that we have
\be\label{defTPn}
\tilde{\tilde{P}}_{n+S+1}\,(x)={Q(x)}W(x)^{-1}\, \frac{d^{m+1}}{dx^{m+1}}\, \left[\, x^i\, Q(x)^m\, W(x)\,\right]
\ee

Using the integrating factor $W$ equation (\ref{corprototype}) can be rewritten as
\be\label{reeq}
\frac{d}{dx}\left[\, W(x)\bfy(x)\,\right]\, =\, \tilde{\tilde{W}}(x)\left[\bfg(x)-\bfphi(x)\right]
\ee
In view of (\ref{defPn}), (\ref{defTPn}) and (\ref{reeq}) we see that if 
$\bfg=\tilde{\tilde{P}}_{n+S+1}$ then $\bfy=P_n$ is a solution of (\ref{prototype}). 

For general polynomial $\bfg$ we only need to expand $\bfg$ in the special polynomials $\tilde{\tilde{P}}_{n}$, as guaranteed by {Proposition}\,\ref{polyno} (note that (\ref{kpBinvert}), (\ref{ttbi}) imply that $\tilde{\tilde{B}}_\infty$ satisfies (\ref{wib})): 
$$\bfg(x)=\sum_{n=0}^{\mbox{deg\ }\bfg}\, \tilde{\tilde{P}}_n(x)\,\bfg_n$$
Letting
$$\bfy(x)=\sum_{n=S+1}^{\mbox{deg\ }\bfg}\, {P}_{n-S-1}(x)\,\bfg_n\ \ \ \ \ \ \ \ {\mbox{and\ \ \ \ \ \ \ \ }}\bfphi(x)=\sum_{n=0}^{S}\tilde{\tilde{P}}_n(x)\,\,\bfg_n$$
the polynomial $\bfy(x)$ satisfies (\ref{corprototype}). \qed

{\bf{Remark.}} Note that deg$\,\bfy(x)=$deg$\, \bfg(x)\,-S-1$.

\subsection{Analytic structure: proof of Lemma\,\ref{FundLem} when the eigenvalues of all matrices $B_j$ have positive real parts.}\label{Step_II}

In \S\ref{locsol} the unique solution of (\ref{prototype}) which is analytic at $x=p_0$ was constructed via its Taylor series. {Lemma}\,\ref{olem} shows that if the eigenvalues of all matrices $B_j$ have positive real parts this solution can be written in an integral form (\ref{anatp0}).

Under the assumptions of \S\ref{Step_II} we can choose the starting point of the integral in (\ref{gsnsolln}) to be a singularity of the equation, and the general solution of (\ref{prototype}) has the form
\be\label{ggg}
W(x)^{-1}\boldsymbol{\xi}+\bfy(x)
\ee
where
\be\label{anatp0}
\bfy(x)\, =\, W(x)^{-1}\,\int_{p_0}^x\, Q(t)^{-1}W(t)\bfg(t)\, dt
 \ee
and $\boldsymbol{\xi}\in\CC^N$.

\begin{Lemma}\label{olem}
If the eigenvalues of $B_0$ have positive real parts, then (\ref{anatp0}) is the unique solution of (\ref{prototype}) which is analytic at $x=p_0$. 
\end{Lemma}

{\em{Proof.}}

Let us first convince ourselves that (\ref{anatp0}) is analytic at $x=p_0$.

There exists a fundamental matrix for (\ref{eqv}) of the form $W_0(x)=(x-p_0)^{\hat{B}_0}\Phi_0(x)$ where $\Phi_0(x)$ is analytic at $p_0$, $\Phi_0(p_0)=I$, and ${\hat{B}_0}$ is a constant matrix (see \S\ref{App1s} for details). Therefore $W(x)=T_0W_0(x)$ for some constant, invertible matrix $T_0$.
Note that (\ref{anatp0}) is the same if we replace $W(x)$ by $W_0(x)$.

Expanding in a convergent series, 
$$Q(t)^{-1}\Phi_0(t)\bfg(t)=\sum_{n=-1}^\infty(t-p_0)^n\bff_n$$
 we see that it is enough to show that 
$(x-p_0)^{-{\hat{B}_0}}\,\int_{p_0}^x\, (t-p_0)^{{\hat{B}_0}+n}\, dt $ is analytic at $p_0$ for all $n\geq -1$. This is clearly true: the usual integration formula for powers applies because
 the matrix ${\hat{B}_0}+n+1$ is invertible due to the assumption (\ref{kpBinvert}) (and any eigenvalue of ${\hat{B}_0}$ is also an eigenvalue of $B_0$).

Uniqueness of the analytic at $p_0$ solution of (\ref{prototype}) was proved in \S\ref{locsol}. Alternatively, a direct proof is by noting that $W(x)^{-1}\boldsymbol{\xi}$ is analytic at $p_0$ only for the constant vector $\boldsymbol{\xi}=0$.
 Indeed, assuming otherwise, $W(x)^{-1}\boldsymbol{\xi}$ is a nonzero solution of $\bfy'+B(x)\bfy=0$ which is analytic at $x=p_0$ and this implies that $-B_0$ has an eigenvalue in $\NN$, which contradicts assumption (\ref{kpBinvert}). \qed

Note that in Lemma\,\ref{olem} we can replace $p_0,B_0$ by any $p_j,B_j$, and these results 
are now used to determine the condition for the solution (\ref{anatp0}) to be analytic at other points $p_j$ as well. Rewriting (\ref{anatp0}) as
\begin{multline}
\bfy(x)
=W(x)^{-1}\int_{p_0}^{p_j}Q(t)^{-1}W(t)\bfg(t)\, dt +W(x)^{-1}\int_{p_j}^xQ(t)^{-1}W(t)\bfg(t)\, dt 
\end{multline}
and using Lemma\,\ref{olem} it follows that $\bfy(x)$ is analytic at the point $p_j$ if only if
\be\label{an_cond}
\int_{p_0}^{p_j}\, Q(t)^{-1}W(t)\bfg(t)\, dt \, =\, 0
\ee

Then $\bfy(x)$ is analytic on $D$ if (\ref{an_cond}) holds for all $j=1,2,\ldots ,S+1$.

Thus a correction $\boldsymbol{\phi}(x)$ needs to be determined so that
\be\label{cond_phi}
\int_{p_0}^{p_j}\, Q(t)^{-1}W(t)\left[\bfg(t)-\boldsymbol{\phi}(t)\right]\, dt \, =\, 0,\ \ \ \ \forall j=1,2,\ldots ,S+1
\ee

Denote $\boldsymbol{\phi}(x)=\sum_{i=0}^{S}\boldsymbol{\phi}_ix^i$. Conditions 
(\ref{cond_phi}) become
$$\sum_{i=0}^{S} \int_{p_0}^{p_j} x^iQ(x)^{-1}W(x)dx \ \boldsymbol{\phi}_i=\int_{p_0}^{p_j}Q(x)^{-1}W(x) \bfg(x)dx,\ j=1,\ldots ,S+1$$

It is enough to show that the system of equation
\be\label{moments}
\sum_{i=0}^{S}\, \mathcal{M}_{ji}\, \, \boldsymbol{\phi}_i=\, \boldsymbol{\xi}_j,\ \ \ j=1,2,\ldots ,S+1
\ee
where
$$\mathcal{M}_{ji}=\int_{p_0}^{p_j}\, x^iQ(x)^{-1}W(x)\, dx \ \ \ {\mbox{and}} \ \ \ \ \ \ \boldsymbol{\xi}_j\in\CC^N$$
has a unique solution $\boldsymbol{\phi}_0,\ldots,\boldsymbol{\phi}_S$.

The linear system (\ref{moments}) has $N(S+1)$ equations and $N(S+1)$ unknowns, and by 
Lemma\,\ref{UniLem} the only solution of the homogeneous system is the zero solution. Therefore the system (\ref{moments}) has a unique solution, completing the proof of Lemma\,\ref{FundLem} under the assumptions of \S\ref{Step_II}.

\subsection{Proof of Lemma\,\ref{FundLem} in the general case: analytic and algebraic methods complementing each other.}
\label{Step_III}

  Lemma\,\ref{Linc} shows that there is an affine transformation of the dependent variable mapping
   (\ref{prototype}) into a similar system which has the eigenvalues of the matrices $B_j$ increased by one:

\begin{Lemma}\label{Linc}

Consider a nonhomogeneous Fuchsian equation (\ref{prototype}). There exists a polynomial $\bfy_0(x)$ of degree at most $S+1$ so that the substitution
\be\label{suby}
\bfy=\bfy_0(x)+Q(x)\tilde{\bfy}
\ee
takes (\ref{prototype}) into the system
\be\label{p1_prototype}
\tilde{\bfy}'(x)+\, \tilde{B}(x)\, \tilde{\bfy}(x)=\, \frac{1}{Q(x)}\, \tilde{\bfg}(x)
\ee
with
\be\label{tildeB}
\tilde{B}(x)=B(x)+\frac{Q'(x)}{Q(x)}\, I,\ \ \ {\mbox{therefore}}\ \tilde{B}_j=B_j+I
\ee
and $\tilde{\bfg}(x)$ analytic on $D$.

\end{Lemma}

{\em{Proof of {Lemma}\,\ref{Linc}.}}

Decompose $\bfg$ as 
\be\label{gpgt}
\bfg(x)=\bfg_P(x)+Q(x)\tilde{\tilde{\bfg}}(x)
\ee
where $\bfg_P$ is a polynomial degree at most $S+1$ and $\tilde{\tilde{\bfg}}$ is analytic on $D$.

Substituting (\ref{suby}) in (\ref{prototype}) we obtain (\ref{p1_prototype}) with
$$\tilde{\bfg}(x)=\tilde{\tilde{\bfg}}(x)-\bfy_0'(x)+\frac{1}{Q(x)}\, \bfg_P(x)-{B}(x)\,\bfy_0(x)$$

The polynomial $\bfy_0(x)$ is now determined so that the difference $\frac{1}{Q(x)}\, \bfg_P(x)-{B}(x)\,\bfy_0(x)$ has no poles, therefore it is a polynomial in $x$, degree at most $S$. Indeed, if $\bfg_j$ are the residues of $\bfg_P/Q$ at $x=p_j$:
\be\label{S1}
\bfg_p(x)=Q(x)\left(\sum_{j=0}^{S+1}\frac{1}{x-p_j}\, \bfg_j\, \right)
\ee
then take
\be\label{S2}
\bfy_0(x)=Q(x)\left(\, \sum_{j=0}^{S+1}\frac{1}{x-p_j}\, \frac{1}{Q'(p_j)}\, B_j^{-1}\,  \bfg_j \right)
\ee
giving (\ref{suby}) and completing the proof of {Lemma}\,\ref{Linc}. \qed

\begin{Remark}\label{Rem}
Note that $\bfy_0$ depends linearly on $\bfg_P$, and it does not depend on $\tilde{\tilde{\bfg}}$.
\end{Remark}

By repeated use of {Lemma}\,\ref{Linc} we can link solutions of systems (\ref{prototype}) with matrices $B_j$ with non-positive real parts to solutions of systems with positive real parts. For such systems the existence of corrections was proved in \S\ref{Step_II}. What is still needed is to show that these corrections can be carried back to the initial system. The following Lemma shows that this is indeed possible. 

\begin{Lemma}\label{PLinc}

Consider an equation (\ref{prototype}) and its transformation (\ref{p1_prototype}) provided by {Lemma}\,\ref{Linc}.

Let $\tilde{\bfphi}$ be a polynomial of degree at most $S$ and consider the following modification of (\ref{p1_prototype}):
\be\label{cp1_prototype}
\tilde{\bfy}'(x)+\, \tilde{B}(x)\, \tilde{\bfy}(x)=\, \frac{1}{Q(x)}\, \left[\tilde{\bfg}(x)-\tilde{\bfphi}(x)\right]
\ee

Then for any $\tilde{\bfphi}$ there exists ${\bfphi}$ a polynomial of degree at most $S$, so that equation 
(\ref{corprototype}) transforms to (\ref{cp1_prototype}) using an affine substitution
\be\label{tsuby}
\bfy=\bfy_{\bfphi}(x)+Q(x)\tilde{\bfy}
\ee
with $\bfy_{\bfphi}(x)$ a polynomial of degree at most $S+1$.

\end{Lemma}

{\em{Proof of {Lemma}\,\ref{PLinc}.}}

It is easy to see (by Remark\,\ref{Rem}, or by direct calculation) that the substitution (\ref{tsuby}) transforms (\ref{corprototype}) to (\ref{cp1_prototype}) if $\bfy_{\bfphi}=\bfy_0-\bfy_1$
with $\bfy_0(x)$ as in the proof of {Lemma}\,\ref{Linc} (see (\ref{gpgt}), (\ref{S1}), (\ref{S2})), and 
$\bfy_1(x)$ a polynomial solution, degree at most $S+1$, of 
\be\label{eqy1}
\bfy'_1+B(x)\bfy_1=\frac{1}{Q}\left[Q\tilde{\bfphi}-\bfphi\right]
\ee

The existence of $\bfphi$ and of the polynomial solution of (\ref{eqy1})  were proved in \S\ref{Step_I}. \qed

\

{\bf{Reduction of the general case to the case of \S\ref{Step_II}.}} 

Let $-n$ (with $n\in\NN$) be a lower bound for the real parts of the eigenvalues of all the matrices $B_0,\ldots,B_{S+1}$.

If $n=0$ this means that \S\ref{Step_II} applies to the system.

If $n=1$ then using {Lemma}\,\ref{Linc} we obtain a system to which \S\ref{Step_II} applies, a correction
$\tilde{\bfphi}$ is found, and using {Lemma}\,\ref{PLinc} a correction $\bfphi$ is obtained for the original system.

If $n\geq 2$ then one uses {Lemma}\,\ref{Linc} repeatedly, $n$ times, then \S\ref{Step_II} applies, a correction is found, and using {Lemma}\,\ref{PLinc} repeatedly, $n$ times, a correction is found for the original system.

The proof of Lemma\,\ref{FundLem} is now complete.

\section{Appendices}

Sections \S\ref{sleq}-\S\ref{App1s} summarize a few known facts regarding solutions of linear, first order differential equations, near an isolated singularity  of the equation. For proofs see \cite{Coddington_Levinson}, followed here except for notation.

\subsection{Singularities of linear differential equation}\label{sleq}

Consider a linear system 
\be\label{ling}
\frac{d\mathbf{y}}{dx}=M(x)\, \mathbf{y},\ \ \ \bfy\in\CC^d
\ee
If the matrix $M(x)$ is analytic at $x_0$ then this point is called a regular point of the equation (\ref{ling}). It is well known that any first order systems (linear or nonlinear) are analytically equivalent to each other near a regular point.

 If the matrix $M(x)$ is not analytic at $x_0$ then $x_0$ is called a singular point of the equation.

 If the singular point $x_0$ of $M(x)$ is isolated\footnote{By an isolated singularity it is meant here that $M(x)$ is analytic in a punctured disk $0<|x-x_0|<\delta$. In particular $M(x)$ is single-valued at $x_0$.} then every fundamental matrix $Y(x)$ of (\ref{ling}) has the form 
 \be\label{fdsoliso}
 Y(x)=\Phi(x)(x-x_0)^P
 \ee
  where $\Phi(x)$ has an isolated singularity at $x_0$ and $P$ is a constant matrix.
 
 {\em{Classification based on the regularity of solutions.}}
 If $\Phi(x)$ has at most a pole singularity at $x_0$, then the point $x_0$ is called a {{regular singularity}}. Otherwise it is called an {{irregular singularity}}. 
 
 Note that if $x_0$ is a regular singularity then, from (\ref{fdsoliso}), the fundamental matrix $Y(x)$ has a convergent expansion as a sum of series in integer powers of $x-x_0$ possibly multiplied by noninteger powers of $x-x_0$ and integer powers of $\ln(x-x_0)$.
 
{\em{Classification based on the regularity of the equation.}} Assume the isolated singularity $x_0$ of $M(x)$ is a pole: $M(x)=(x-x_0)^{-r-1}L(x)$ with $L(x)$ is analytic at $x=x_0$ and $r\geq 0$. The number $r$ is called the rank of the singularity.

If $r=0$ the point $x_0$ is called a singular point of the first kind, or a Fuchsian point. It turns out that a Fuchsian point is necessarily a regular singuarity. 

More details about the structure of solutions at a Fuchsian point are discussed in the following \S\ref{App1s}.

\subsection{Systems near a Fuchsian point}\label{App1s}

Consider a linear system 
\be\label{lin_one_sing}
\frac{d\mathbf{y}}{dx}=\frac{1}{x-x_0}L(x)\, \mathbf{y},\ \ \ \bfy\in\CC^d
\ee
with $L(x)$ a matrix analytic at $x_0$. 

The structure of the matrix $P$ in (\ref{fdsoliso}) is particularly simple if the eigenvalues of $L(x_0)$ are linearly {{nonresonant}}, in the sense that 
\be\label{lnr}
{\mbox{any\ two\ {{different}}\ eigenvalues\ do\ not\ differ\ by\ an\ integer}}
\ee
In this case we can take $P=L(x_0)$ and choose a fundamental matrix in the form
\be\label{genforsol}
Y(x)= \Phi(x)\,(x-x_0)^{L(x_0)}\ \ \ \ {\mbox{with}}\ \Phi(x_0)=I
\ee
and $\Phi(x)$ analytic at $x_0$.

In linearly resonant cases (when there are distinct eigenvalues of $L(x_0)$ which differ by an integer) the matrix $P$ may differ from $L(x_0)$. However, the eigenvalues of $P$ (not counting their multiplicity) form a subset of the eigenvalues of $L(x_0)$\footnote{Any resonant string of eigenvalues of $L(x_0)$ is replaced in $P$ by the lowest eigenvalue in the string with higher multiplicity, see in \cite{Coddington_Levinson} the Lemma of Chap.4, Sec.4.}.

Another representation of the general solution in resonant cases is again (\ref{genforsol}), only that this time $\Phi(x)$ is a convergent power series in both $x-x_0$ and $\ln(x-x_0)$.

Note that in (\ref{eqv}) there is a different order of multiplication in the differential equation, and the Fuchsian point in (\ref{chooseW}) is $x_0=\infty$.

\subsection{{Nonlinear perturbation of first order systems with a nonresonant Fuchsian 
point.}}\label{nl1fuchs}

Consider a nonlinear equation
\be\label{non_one_sing}
\frac{d\mathbf{u}}{dx}=\frac{1}{x-x_0}L(x)\, \mathbf{u}+\frac{1}{x-x_0}\bff(x,\bfu),\ \ \ \bfu\in\CC^d
\ee
with $\bff(x,\bfu)$ analytic at $x=x_0$, with a zero of order two at $\bfu=0$ and $L(x_0)$ linearly nonresonant.

Let $\Phi(x)$ be as in (\ref{genforsol}). The substitution $\mathbf{u}=\Phi(x)\tilde{\mathbf{u}}$ in (\ref{non_one_sing}) yields
\be
\frac{d\tilde{\mathbf{u}}}{dx}=\frac{1}{x-x_0}L(x_0)\tilde{\mathbf{u}}+\frac{1}{x-x_0}\Phi(x)^{-1}\bff(x,\Phi(x)\tilde{\mathbf{u}})
\ee
which is an equation with the same regularity as (\ref{non_one_sing}) for small $x-x_0$ and $\mathbf{u}$ and with a constant matrix in the linear part.

\subsection{Proof of Lemma\,\ref{sollem}.}\label{proofLJ}

This result appears in \cite{Arnold} Ch.5, \S{22}, where it is proved for $M$ diagonal, and it is stated for $M$ in Jordan normal form. A proof is included here for completeness.

{\bf{I.}} Note that  (after a change of coordinates in $\CC^d$) we can assume that $M$ is in normal form. 

Indeed, upon a change of coordinates in $\CC^d$, taking $\bfw$ to $S\bfw$, a $\CC^d$-valued polynomial $\bfq(\bfw)$, $\bfq\in\mathcal{P}_n$ becomes $S^{-1}\bfq(S\bfw)$. Denote by
 ${S}_\#$ this induced change of coordinates on $\mathcal{P}_n$: the operator $J_M$ becomes, in the new coordinates, ${S}_\# ^{-1}J_M{S}_\#$, which equals $J_{S^{-1}MS}$ (as seen after a short calculation).

{\bf{II.}} Note that if $M$ invariates the span of $\{{\rm{e}}_k;s\leq k\leq t\}$ (for some $1\leq s\leq t\leq d$) then $J_M$ invariates the span of $\{\mathbf{r}_{\bfm,k};s\leq k\leq t,|\bfm|=n\}$ - see (\ref{canbasP}) for notation.

This follows because $({\rm{d}}_\bfw\mathbf{r}_{\bfm,k})M\bfw$ is a (polynomial) multiple of ${\rm{e}}_k$:
\be\label{pmej}
({\rm{d}}_\bfw\mathbf{r}_{\bfm,k})M\bfw=\phi_\bfm(\bfw){\rm{e}}_k
\ee


{\bf{III.}} Assume $M$ is in normal Jordan form. This means that there is a partition $1=t_0<t_1<\ldots<t_p=d$ so that 
\be\label{Mej}
M{\rm{e}}_j=\lambda_j{\rm{e}}_j+\sigma_j{\rm{e}}_{j-1}
\ee
 with $\lambda_j$ equal for all $j$ with $t_i\leq j<t_{i+1}$, and with $\sigma_{t_i}=0$, and all other $\sigma_j=1$.

Using (\ref{Mej}) in (\ref{pmej}) we obtain the formula for the polynomial $\phi$:
\be\label{fp}
\phi_\bfm(\bfw)=(\boldsymbol{\lambda}\cdot\bfm)\bfw^\bfm+\sum_{j=1}^{d-1}\sigma_{j+1}m_j\bfw^{\bfm-{\rm{e}}_j+{\rm{e}}_{j+1}}
\ee

{\bf{IV.}} The remarks above imply that the matrix of $J_M$ in the canonical base (\ref{canbasP}) is block diagonal. Moreover, with an appropriate ordering of this base the matrix of $J_M$ is upper diagonal.

This ordering is as follows: we have $(\bfm',k')\preccurlyeq(\bfm,k)$  if $k'<k$, or $k'=k$ and $\bfm'\leq\bfm$ in lexicographic order.

As a example, take $d=3$. For $k=1$, we have in increasing order\footnote{As a reminder, only $\bfm\in\NN^d$ with $|\bfm|=n$ are present.}: $(0,0,n)$, $(0,1,n-1)$, $(0,2,n-2)$,\ldots , $(0,n,0)$, $(1,0,n-1)$, $(1,1,n-2)$, $\ldots\ (n,0,0)$, after which comes the same sequence with $k=2$ etc.

The fact that the matrix of $J_M$ is upper diagonal is expressed by the fact that $J_M\mathbf{r}_{\bfm,k}$ is a linear combination of $\mathbf{r}_{\bfm',k'}$ with $(\bfm',k')\preccurlyeq(\bfm,k)$.

To show this, consider the first block for $M$: $k$ from $1$ to $t_1$. 

For $k=1$ the first vector of the canonical basis is $\mathbf{r}_{(0,\ldots,0,n),1}$. From (\ref{fp}) we have $\phi_{(0,\ldots,0,n)}(\bfw)=(\boldsymbol{\lambda}\cdot\bfm)\bfw^\bfm$ and from (\ref{Mej}) it follows that $M\mathbf{r}_{(0,\ldots,0,n),1}=\lambda_1\bfw^\bfm{\rm{e}}_1$ hence $\mathbf{r}_{(0,\ldots,0,n),1}$ is an eigenvector of $J_M$ corresponding to the eigenvalue $\boldsymbol{\lambda}\cdot\bfm-\lambda_1$. For other values of $\bfm$ it follows from (\ref{pmej}), (\ref{Mej}), (\ref{fp}) that $J_M\mathbf{r}_{\bfm,1}$ is a linear combination of the vectors $\mathbf{r}_{\bfm,1}$ (with coefficient $\boldsymbol{\lambda}\cdot\bfm-\lambda_1$) and $\mathbf{r}_{\bfm-{\rm{e}}_j+{\rm{e}}_{j+1},1}$, $j=1,\ldots,d-1$. We obviously have $(\bfm-{\rm{e}}_j+{\rm{e}}_{j+1},1)\preccurlyeq (\bfm,1)$.
(Note that the index ${\bfm-{\rm{e}}_j+{\rm{e}}_{j+1}}$ is not defined if $m_j=0$, or if $m_{j+1}=n$. Since in both cases $m_j=0$ these indexes do not appear in (\ref{fp}).) 

For $k=2,\ldots,t_1$ the discussion is similar, only, due to (\ref{Mej}) also earlier terms, with $k=1$, appear.  

For the other blocks of $M$ the proof is similar.\qed

\end{document}